\documentclass[11pt]{article}
\usepackage{amsmath}
\usepackage{graphicx}
\usepackage{amsfonts}
\usepackage{amssymb}
\usepackage{mathrsfs}

\newtheorem{theorem}{Theorem}
\newtheorem{lemma}[theorem]{Lemma}

\newtheorem{definition}[theorem]{Definition}

\newcommand\beq{\begin{equation}}
\newcommand\eeq{\end{equation}}
\newcommand\bqu{\begin{quote}}
\newcommand\equ{\end{quote}}

\def\C{\mathbb C}

\def\ra{\rightarrow}

\def\neq{\not=}

\def\inv{^{-1}}

\def\jump{{\hspace{1cm}}}

\def\la{\lambda}

\def\ga{\gamma}

\def\ka{\kappa}

\def\BM{\begin{matrix}}
\def\EM{\end{matrix}}
\def\BpM{\begin{pmatrix}}
\def\EpM{\end{pmatrix}}

\def\squarebox#1{\hbox to #1{\hfill\vbox to #1{\vfill}}}
\newcommand{\qed}{\hspace*{\fill}
\vbox{\hrule\hbox{\vrule\squarebox{.667em}\vrule}\hrule}\smallskip}

\begin{document}
\bibliographystyle{plain}
\renewcommand{\thesection}{\Roman{section}}
\thispagestyle{empty}

\begin{center}
{\Large\bf Spectral data and solvent theory

for regular matrix polynomials}\\[0.8cm]
{\bf Nir Cohen, Edgar Pereira}\footnote{Department Mathematics,
UFRN - Universidade Federal de Rio Grande do Norte, Natal, Brazil
(nir@ccet.ufrn.br, edgar@ccet.ufrn.br).}
\\[0.6cm]
\today\\[1.0cm]
\end{center}

\bqu\noindent {\bf Keywords:~~} {\em matrix polynomial, bisolvent,
solvent, standard pair, Kronecker-Weierstrass form, strict
equivalence, Jordan chain, right linear factors}

\bqu \noindent {\bf Abstract:~~} {\em This paper contains a
re-evaluation of the spectral approach and factorizability for
regular matrix polynomials. In addition, solvent theory is
extended from the monic and comonic cases to the regular case. The
classification of extended solvents (bisolvents) is shown to be
equivalent to the classification of all the regular first order
right factors for a general matrix polynomial.

}\equ

\section{Introduction}

It is of interest to extend the fundamental theorem of algebra
(roots and factorization) from complex polynomials to complex
matrix polynomials, i.e. expressions of the form
$$P(\la)=\sum_{i=0}^k A_i\la^i\,\jump A_i\in M_{n,m}(\C),\quad\la\in\C.$$

{\bf (I) KW theory.} This 19th century approach, which is limited
to matrix pencils, extends the usual notion of complex
eigenvalues, including roots "at infinity". Initially, Weierstrass
considered only {\em regular} pencils, and provided a complete
analysis of roots (also called eigenvalues) via strict equivalence
and the Weierstrass (or W) canonical form. Kronecker extended
Weierstrass' work to {\em singular} pencils in terms of additional
structure: right and left null indices (\cite{gan}, p. 25) , later
described in detail by Forney \cite{for}. These considerations
have led to the Kronecker-Weierstrass (KW) form under strict
equivalence (\cite{gan}, p. 24).

{\bf (II) Standard pair theory.} Although complex roots can be
defined for higher order matrix polynomials, the associated
analysis of "eigenvectors" in terms of state space models became
possible only after the emergence of the state space approach in
engineering in 1960, due to Kalman \cite{kal}. A series of papers
in the 1970-80s provided a complete state space analysis for
regular matrix polynomials in terms of strong equivalence and
reduction to a companion pencil. The eigenvector data and the
eigenvalue data are stored in two separate matrices
$([X,Y],T\oplus Z)$
 called standard pair \cite{nir} or decomposable pair
\cite{goh}. Notably, these matrices satisfy the following algebraic relations:
 $$\sum_{i=0}^k A_iXT^i=0,\jump \sum_{i=0}^k A_iYZ^{k-i}=0$$
plus a certain controllability condition. The matrix polynomial
can be reconstructed in an essentially unique way from this data.
See \cite{rod-mon}, \cite{rod-com}, \cite{rod-com2}, \cite{nir},
\cite{goh}, \cite{ler}, \cite{lanc-res}.

{\bf (III) Factorization.} The description of factorizations of
the type $P(\la)=Q(\la)R(\la)$ are far more difficult, mainly due
to the commutativity issue. Only the monic $2\times 2$ case is
fully analyzed \cite{nir2x2}. Some results are known for finding a
first order right factor, specially in the case of quadratic
matrix polynomials, see \cite{lanc-div}, \cite{goh}.

{\bf (IV) Solvent theory.} Solvent theory has a different, truly
non-commutative starting point, where roots are not complex numbers but $m\times m$ matrices. Formally, considering $P(\cdot)$ as a left multiplication operator on $M_m(\C)$, solvents are defined as solutions $S\in M_{m}(\C)$ of the matrix equation
 $$\sum_{i=0}^k A_iS^i=0.$$
The scope of this theory is, to date, quite limited: the
classification of solvents is complete only in the essentially
monic case $(n=m$ and $det~A_k\neq 0)$ \cite{per}. Nevertheless,
solvent theory shows close ties with both the standard pair
approach and factorization, especially first order right factors.


%
%


A similar classification of "cosolvents" is available in the
comonic case, based on the "reverse equation"
 $$\sum_{i=0}^k A_iS^{k-i}=0,$$
leading to a mirror image of cosolvent theory.

%
%
%

The objective of this paper is to integrate the solvent and cosolvent
approaches in a non-trivial way to a unique theory which is able to describe all the
{\em regular} first order right factors of an arbitrary matrix polynomial.
We use the concept of {\em bisolvent}, related with the algebraic relation
 $$\sum_{i=0}^k A_iS_1^iS_2^{k-i}=0.$$
Clearly, commutativity is a subtle issue here, and so we find it
necessary to restrict our scope to commuting pairs $(S_1,S_2)$; as
we show below, however, there is absolutely no need to consider
non-commuting bisolvents!

In fact, we focus in this paper a special type of commuting
bisolvents, which we call "separable". These are, formally,
triples $(S_1,S_2,\Pi)$ with $\Pi\in M_m(\C)$ idempotent, which
also satisfy
 \beq\label{idem}S_1=\Pi S_1\Pi+(I-\Pi),\jump S_2=\Pi+(I-\Pi)S_2(I-\Pi).\eeq

The text is self sufficient, but the first few preliminary
sections are quite brief as they mostly summarize existing
results.


\section{Strict equivalence and KW theory}

The following basic definitions will be used throughout.

\begin{quote}\begin{definition}\label{basic}
\quad (i) $M_{nm}$, $M_n$ and $GL_n$ will denote the sets of
$n\times m,$ $n\times n$ and non-singular $n\times n$ (complex)
matrices.

\quad (ii) A (complex) $n,m,k$-polynomial is an expression of the form $P(\la)=\sum _{i=1}^kA_i\la^i$ with $A_i\in M_{nm}$ and $A_0,A_k\neq 0.$

\quad (iii) $P(\la)$ is called monic if $A_k=I_n$, comonic if $A_0=I_n,$ regular if $n=m$ and $det~P(\la)\neq 0$ (singular otherwise), unimodular if $det~P(\la)$ is a non-zero constant, i.e. having a polynomial inverse.

\quad (iv) The reverse $n,m,k$ polynomial is $\tilde P(\la)=\la^kP(1/\la)=\sum_{i=0}^k A_i\la^{k-i}.$
\end{definition}\end{quote}

Although in our discussion we include considerations on singular
matrix polynomials the scope of the presented theory here in
general is limited to regular matrix polynomials.

A central question is the classification of matrix pencils up to
strict equivalence.

\begin{quote}\begin{definition}\label{strict}
\quad (i) Two $n,m,k$ polynomials $P(\la)$ and $P'(\la)$ are called strictly equivalent if $QP(\la)=P'(\la)R$ with $Q\in GL_n$ and $R\in GL_m.$

\quad (ii) An $n,m,1$ polynomial is called an $n,m$ pencil.

\quad (iii) A Weierstrass pencil is a direct sum of monic and comonic Jordan pencils.

\quad (iv) Denote by $R_\ka(\la)$ the matrix pencil formed from
the first $\ka$ rows of the pencil $\la I+J$, where $J$ is a
single nilpotent Jordan block of size $(\ka+1).$ We call $R_\ka$
(resp. $R^T_\ka(\la))$ a right null block (resp. left null block)
of index $\ka$.

\quad (v) A KW pencil is a direct sum of the form
 $$T(\la)=\BpM T'(\la)&0\\ 0&0\EpM,\jump T'(\la)=\oplus T_i(\la)$$
consisting of four types of units $T_i(\la)$: monic and comonic Jordan pencils,
right and left null blocks.
\end{definition}\end{quote}

\noindent The following is the main result of the KW theory.

\begin{quote}\begin{theorem}\label{sing}
\quad (i) Every matrix pencil $C(\la)$ is strictly equivalent to a KW pencil $T(\la).$

(ii) $C(\la)$ is regular iff $T(\la)$ is Weierstrass. In addition, $P(\la)$ is essentially monic,
 i.e. $det~A_k\neq 0$ (resp. essentially comonic) iff $T(\la)$ is monic (resp. comonic).

(iii) The KW form is unique up to two trivial operations: (i) reordering the blocks;
(ii) replacing a monic Jordan block of eigenvalue $\la\neq 0$
by a comonic Jordan block of the same size for $1/\la,$ and vice versa.
\end{theorem}\end{quote}

To avoid the last source of non-uniqueness, it is commonplace to require that all
comonic Jordan blocks be nilpotent; namely, the monic part includes all the
finite eigenvalues.

The KW form defines the full set of strict equivalence invariants: the eigenvalues
(eigenvalues of monic blocks plus inverse eigenvalues of comonic blocks), their partial multiplicities, right and left null indices.

\section{Strong equivalence and linearization}

We now discuss the extension of the KW theory to higher order polynomials.

\begin{quote}\begin{definition}\label{strict}
(i) We say that $P(\la)$ and $P'(\la)$ (of possibly different sizes $n\times m$ and $n'\times m')$) are equivalent if
there exist unimodular matrix polynomials $Q(\la),R(\la)$ and integers $p,q$ so that
 $$[P(\la)\oplus I_p]Q(\la)=R(\la)[C(\la)\oplus I_q].$$
Clearly, $n-m=n'-m'.$

(ii) $P(\la)$ and $P'(\la)$ are said to be strongly equivalent $(P\sim P')$
if both the direct pair $P,P'$ and reverse pair $\tilde P,\tilde P'$ are equivalent.

(iii) $C(\la)$ is a linearization for $P(\la)$ if $C(\la)$ is a matrix pencil and $C\sim P$.
\end{definition}\end{quote}


\begin{quote}\begin{theorem}\label{facts-pencil}

\quad (i)\quad  Strict equivalence implies strong equivalence.

\quad (ii) \quad Two matrix pencils are strongly equivalent iff they have
the same eigenvalues and partial multiplicities and the same number (but not necessarily values)
 of left and right null indices.

\quad (iii) So, for regular pencils strong equivalence implies strict equivalence.
\end{theorem}\end{quote}

Theorem \ref{facts-pencil}(ii) does not extend to singular pencils.
In fact, strongly equivalent singular pencils may have different null indices.

\quad Order reduction and companion forms are used in extending
Theorem \ref{facts-pencil} to higher order. Coefficient-based
order reduction is a powerful method in ODEs and engineering
applications. For monic matrix polynomials, four monic companion
forms are widely used in the engineering literature, which may be
described as "up, down, right" and "left" \cite{kai}, \cite{bar};
throughout this work we restrict attention to the {\em down and
right companion pencils} defined by

$$C_d(\lambda;P)={\tiny\BpM
I      & 0 & \cdots &       0 &      0   \\
0      & I & \cdots &       0 &      0   \\
\vdots &   &        &         & \vdots   \\
0      & 0 & \cdots &       I &      0   \\
0      & 0 & \cdots &       0 &      A_m \\
\EpM}\lambda - {\tiny\BpM
0      & I   & \cdots &       0 &      0       \\
0      & 0   & \cdots &       0 &      0       \\
\vdots &     &        &         & \vdots       \\
0      & 0   & \cdots &       0 &      I       \\
-A_0    & -A_1 & \cdots & -A_{m-2} &      -A_{m-1}\EpM},$$
$$C_r(\la;P):=C_d(\la;P^T)^T$$ in which the coefficients $A_i$
occupy, respectively, the lowest block row and the rightmost block
column. In general the sizes of these two pencils are different:
resp, $mk\times [m(k-1)+n]\times mk$ and $[n(k-1)+m]\times nk$.
Other companion forms, besides these four, are discussed in
\cite{fie} and \cite{ant}. A pencil is considered as its own
linearization, i.e. $C_d(\la;P)=C_r(\la;P)=P(\la).$

\begin{quote}\begin{theorem}\label{results-comp}
\quad (i) \quad Any matrix polynomial $P(\la)$ admits linearizations; for example, all the companion forms including $C_d(\la;P)$ and $C_r(\la;P)$.

(ii) If the $n,n,k$ polynomial $P(\la)$ is regular then all its regular linearizations
are of size $nk,$ regular, and pairwise strictly equivalent.

(iii) If $P_i(\la)$ have equal degree and $C_i(\la)$ are linearizations
for $P_i(\la)$ then $\oplus C_i$ is a linearization for $\oplus P_i.$
 \end{theorem}\end{quote}

\noindent Item (iii) may fail for polynomials of different degrees, due to
indiscrepancies at the infinite eigenvalue. For example, a direct sum of
monic polynomials of different degrees is not monic, and so
does not have a monic linearization.

Item (ii) fails in the singular case, as $P(\la)$ admits linearizations of various sizes.
\section{The Jordan chain approach}

In order to encode spectral information
efficiently we need to translate the data in matricial rather than modular terms.
Namely, using Taylor expansion of root vectors at a neighborhood of a root, we obtain certain chains of "eigenvectors" and "generalized eigenvectors". With this,
a direct generalization of the usual concept of Jordan chains of a matrix is
obtained. The problem is that the matricial representation is adapted to encode eigenvector data but not null vector data.

\begin{quote}\begin{definition}\label{chain}
(i) A Jordan $\ka$-chain for $P(\la)$ at $a\in\C$ is a sequence
$v_1,\cdots,v_\mu$ with $v_1\neq 0$ so that
 $$\frac{1}{j!}P^{(j-1)}(a)v_1+\cdots+P(a)v_j=0\jump(0\leq j\leq \ka-1).$$

(ii) A Jordan $\ka$-cochain at $\la=\infty$ for $P(\la)$ is a Jordan
$\ka$-chain for $\tilde P(\la)$ at $\la=0.$

(iii) A Jordan system is a finite set of chains and cochains associated with complex points
in which the eigenvectors (initial vectors) in each chain
associated with $\la=a$ are LI.\end{definition}\end{quote}

\noindent In the regular case, this procedure makes sense and it
is clear that a Jordan system is limited in size. Indeed, the
number of chains at each eigenvalue is limited by the linear
independence requirement, and the degree of each chain is limited
due to size of respective Jordan block.

\begin{quote}\begin{lemma}\label{indices-chain}Let $P(\la)$ be regular.

\quad (i) The chain lengths due to a maximal Jordan system
coincide with the partial multiplicities at $\la$ in the Kronecker form.

(ii) The chain lengths of a non-maximal Jordan system for $P(\la)$
are majorized by these.
\end{lemma}\end{quote}

This will be clarified in the next section.

\section{Matricial spectral analysis}


In the regular case both the companion and KW linearizations are
of the same size, hence they are strictly equivalent. Thus,a
maximal system of Jordan chains can be calculated from the strict
equivalence between them. This procedure, developed in the
1970-80s (see \cite{rod-mon},\cite{rod-com},\cite{nir},\cite{goh})
partly goes back to work of Keldy\v s on ODEs in the 1950s
\cite{kel}, although we find those Taylor's expansions, "Jordan
chains" (see Roth \cite{rot}), in a previous work with solutions
for "unilateral equation in matrices", that is $\sum_{i=0}^k
A_iX^i=0,$ for a rectangular matrix $X$. The following analysis
using the term "standard pair" is due to \cite{nir}, see also
\cite{goh} where the term "decomposable pair" is used instead.

When Jordan chains are encoded matricially, some algebraic
relations become evident.

\begin{quote}\begin{definition}\label{def-eigen}
An standard pair for $P(\la)$ is a pair $([X,Y],T\oplus Z)$ where

(i) $X\in M_{m,p}$, $T\in M_{p,p}$ satisfy the algebraic relation
$\sum_{i=0}^k A_iXT^i=0.$

(ii) $Y\in M_{m,q}$, $Z\in M_{q,q}$ satisfy the algebraic relation
$\sum_{i=0}^k A_i YZ^{k-i}=0.$

(iii) the matrix
 \beq\label{cont} Q_k=
 {\tiny\BpM X            &      YZ^{k-1}   \\
XT         &      YZ^{k-2}   \\
\vdots         &      \vdots                     \\
XT^{k-2} &      YZ           \\
XT^{k-1} &      Y                   \\
\EpM}\in M_{nk,p+q}\eeq has full row rank $(=p+q)$.

(iv) We call the pair Jordan if $T\oplus Z$ is in Jordan form.
\end{definition}\end{quote}

The matrix $Q_k$ is sometimes called controllable matrix, this is
due to the fact that a pair satisfy item (iii) above if and only
the matrix $Q_k$ contructed with any $j>1$ instead of $k$, is of
full rank, such pair is called controllable.




We observe that if the dimension of $Z$ is zero ($q=0$), then
$(X,T)$ is a standard pair for $P(\la)$, conversely if the
dimension of $T$ is zero ( $p=0$), then $(Y,Z)$ is a standard pair
for $P(\la)$. Furthermore if $(X,T)$ is a standard pair for
$P(\la)$ with $T$ being nonsingular then $(X,T^{-1})$ is a
standard pair for $\tilde P(\la),$ equivalent conclusions are
obtained if $Z$ or both $T$ and $Z$ are nonsingular. In the case
that the pair $([X,Y],T\oplus Z)$ is Jordan, with $T\oplus Z$
being nonsingular and $(T\oplus Z)^{-1}=R(\tilde T \oplus \tilde Z
)R^{-1}$, where $\tilde T \oplus \tilde Z$ is also in Jordan from,
then $([\tilde X,\tilde Y], (\tilde T \oplus \tilde Z )$ is a
standard pair for $\tilde P(\la)$, where $\tilde X= XR$ and
$\tilde Y=YR$.

{\bf Example 1.} The special case of a first order monic
polynomial is suggestive. An standard pair $(v,a)$ for an $n\times
n$ matrix $A$ induces a Jordan pair $(X=v$, $J=a)$ for $P(\la)=\la
I-A$. The algebraic relation in question reduces to $Av-Iv\la=0.$
\qed

Similarity for standard pairs is defined, quite naturally, in
terms of the gauge transformation
 $$([X_1,X_2],T_1 \oplus T_2)\ra([X_1G,X_2H],G\inv T_1G \oplus H\inv T_2H)\jump (G\in GL_p,H\in GL_q),$$
and it is clear that the entire similarity orbit of an standard
pair for $P(\la)$ consists of standard pairs. However, there
exists an additional degree of freedom not present in the monic or
comonic cases: spectral inversion.

\begin{definition}\label{sp-inv}Spectral inversion of a standard pair $([X_1,X_2],T_1 \oplus T_2)$ is defined
 as follows: assume that $X_1=[X'_1,X'_2]$ and $T_1=T'_1\oplus T'_2$ with $T'_2$ non-singular. Then
 spectral inversion replaces $([X_1,X_2],T_1 \oplus T_2)$ by $([X'_1,[X'_2,X_2]],T'_1 \oplus (T'^{-1}_2 \oplus T_2)).$
  Namely, $p$ is decreased and $q$ is increased, and spectral content of $(X_1,T_1)$ is transferred
   to $(X_2,T_2).$
\end{definition}

The inverse operation considered before is also a spectral
inversion. The zero eigenvalues of $T$ and $Z$ do not suffer
spectral inversion.

Large standard pairs can be constructed from smaller ones by
merging:
 $$X=[X_1,\cdots,X_t],\quad T=\oplus_{i=1}^t T_i,\quad Y=[Y_1,\cdots,Y_t],
 \quad Z=\oplus_{i=1}^tZ_i.$$
These maximal standard pairs correspond to the maximal Jordan
systems discussed in the last section.



If we require that $Z$ be nilpotent, this means that $T$ will have
all the information about the finite spectra and $Z$ all the
information of the infinite spectra of $P(\la)$.

\begin{quote}\begin{theorem}\label{spect-to-pair}
\quad (i) Every regular $n,k$ $P(\la)$ admits a standard pair of
size $nk$.

(ii) $[X,Y],T\oplus Z)$ defines $P(\la)$ via an equivalence
relation with its companion pencil $C_d(\la;P)$. In particular,
the finite eigenvalues, and respective partial multiplicities, of
$P(\la)$ match with the union of eigenvalues of $T$ and inverse
eigenvalues of $Z$, and the infinite eigenvalue, with respective
partial multiplicities, of $P(\la)$ matches with zero eigenvalue
of $Z$.

\quad (iii) $P(\la)$ determines uniquely $([X,Y],T\oplus Z)$ up to
similarity and spectral inversion (in particular, assuming $Z$
nilpotent, up to similarity).

\quad (iv) $([X,Y],T\oplus Z)$ is a Jordan pair (i.e. $T\oplus Z$
is Jordan) iff $[X,Y]$ defines a complete system of Jordan chains
for $P(\la).$
\end{theorem}\end{quote}

\noindent Indeed, consider the strict equivalence $C_d(\la;P)Q=RT(\la)$
between the companion $C_d(\la;P)$ and any separable
linearization $T(\la)=(\la I-T)\oplus(\la Z-I)$ (e.g. the Weierstrass form).
By direct calculation $R,Q$ must have the form
 \begin{equation}\label{}
Q={\tiny\BpM X            &      YZ^{m-1}   \\
XT         &      YZ^{m-2}   \\
\vdots         &      \vdots                     \\
XT^{m-2} &      YZ           \\
XT^{m-1} &      Y                   \\
\EpM},
\jump R={\tiny\BpM
X            &      YZ^{m-2}   \\
XS         &      YZ^{m-3}   \\
\vdots         &      \vdots                     \\
XS^{m-2} &      Y            \\
A_mXS^{m-1} &     -\Sigma^{m-1}_{i=0}A_iYZ^{m-1-i}
\EpM}\end{equation}
for some data
 $$X\in M_{m,p},\quad Y\in M_{m,q},\quad T\in M_p,\quad Z\in M_q\jump (p,q \in{\mathbb  N}).$$
In particular, if $T,Z$ are Jordan matrices, columns of $X,Y$ form a complete set of Jordan chains for $P(\la).$\qed

Let $\ga$ be a simple closed Jordan curve system in $\C$. If we
seek only standard pairs separated by $\ga$, i.e. eigenvalues of
$T$ are in $int(\ga)$ and inverse eigenvalues of $Z$ are in
$ext(\ga),$ then spectral inversion is not an option and the
uniqueness is up to similarity only. For example, standard pairs
with $Z$ nilpotent.

Clearly, in Theorem \ref{spect-to-pair} $P(\la)$ is essentially monic (resp. essentially comonic)
 iff (up to spectral inversion) $q=0$ (resp. $p=0).$

{\bf Example 2.} Consider the matrix $A$ in Example 1. If $J$ is
the Jordan form of $A$, $P(\la)=\la I-A$ has the monic Weierstrass
pencil $\la I-J$ and the monic companion $\la I-A$. Let $X$ be the
similarity matrix between the two pencils, $\la (I-A)X=X(\la I-
J).$ As is well known, the columns of $X$ provide a complete set
of Jordan chains for $A$. Thus, $(X,J)$ is a standard pair. Here
$X=Q_1$ serves as a controllability matrix. \qed



\section{Solvent theory}

An $m\times m$ matrix is a solvent (resp. cosolvent) for
$P(\la)=\sum_{i=0}^k A_i\la^i$ if it satisfies $\sum_{i=0}^k
A_iS^i=0$ (resp. $\sum_{i=0}^k A_iS^{k-i}=0).$ Solvents and
cosolvents help describe right first order factors of $P(\la)$
which are, respectively, monic and comonic.

The solvents can be obtained from a standard pair $([X,Y],T \oplus
Z)$ when $Z$ is nilpotent, that is, $T$ is similar to the finite
part and $Z$ is similar to the infinite part of the canonical
(Weierstrass) form of $P(\la)$, in such case for the sake of
simplicity we will refer $T$ as being the finite part and $Z$ as
being the infinite part of the standard pair.

\begin{quote}\begin{theorem}\label{monic-solv-thm}
\quad (i)\quad  All the solvents $S$ of $P(\la)$ can be
constructed from a maximal controllable standard pair $([X,Y],T
\oplus Z)$, with $Z$ being nilpotent, by  a restriction to an
$m$-dimensional $T$-invariant subspace. Namely, $S=\tilde X\tilde
T\tilde X\inv$ where $(\tilde X,\tilde T)$ are the restrictions of
$(X,T)$ to the subspace.

\quad (ii) \quad In particular, the spectrum of $S$ is subordinate
to that of $T$.

\quad (iii) \quad $\la I-S$ is a right monic factor iff $S$ is a
solvent.
\end{theorem}\end{quote}

The conditions above can be stated in a generalization of
eigenpairs (see \cite{per}), where the necessary and sufficient
conditions for the existence of a solvent are in terms of the
Jordan chains of $P(\la)$ and the partial multiciplities of the
respective eigenvalues.

The number of solvents will depend on the different pairs $(\tilde
X,\tilde T)$ that can be obtained in the theorem above from
$(X,T)$ or equivalent with the different nonsingular matrices
constructed with the leading vectors of the Jordan chains
associated with the finite eigenvalues . This number can be zero,
a finite number or infinitely many.

{\bf Example 3.} Consider the following polynomial $P(\la)= A_1
\la +A_0$, where $A_1$ is singular and $A_0$ is nonsigular. It is
immediate that $P(\la)$ has no solvents. On the other hand, it is
more complicate to construct a polynomial with no solvents for the
monic case (see \cite{per}).\qed

We recall that if the $deg(det(P(\lambda))=M$, then $M$ is the
number of finite eigenvalues, including multiplicities.
Furthermore, if $M=mn$, $P(\lambda)$ has only finite eigenvalues,
otherwise if $M<mn$, $P(\lambda)$ has infinite eigenvalues,
including the multiplicities.

Thus, the maximum finite number of solvents in the case of
$deg(det(P(\lambda))=M$ will be $\left(\begin{array}{cc}
M \\
n \\
\end{array}\right)$, and it occurs when $P(\lambda)$ has $M$ eigenvectors that are linearly
independent $n$ by $n$.

Infinitely many solvents occur when there exists a solvent $X_1$
with an eigenvalue having the geometric multiplicity less in $X_1$
than in $P(\lambda)$.

A similar complete classification of "cosolvents" is based on the
"reverse equation"
 $$\sum_{i=0}^k A_iS^{k-i},$$
leading to a mirror image of cosolvent theory:

\begin{quote}\begin{theorem}\label{com-solv-thm}
\quad (i) \quad All the cosolvents $S$ of $P(\la)$ can be
constructed from a maximal controllable standard pair $([X,Y],T
\oplus Z)$, with $T$ being nilpotent, by restriction to an
m-dimensional $Z$-invariant subspace, in the form $S=\tilde
Y\tilde Z\tilde Y\inv$.

(ii) In particular, the spectrum of $S$ is subordinate to that of
$Z$.

(iii) $I-\la S$ is a right comonic factor iff $S$ is a cosolvent.
\end{theorem}\end{quote}

The relation between solvents and cosolvents is a direct
consequence from their constructions. So, we have that if $S$ is a
nonsingular solvent of $P(\la)$, then $S^{-1}$ is a cosolvent of
$P(\la).$ Moreover, with the same arguments used before, we can
estimate the number of cosolvents and also of "infinite solvents",
that is nilpotent solvents of $\tilde{P}(\la)$.

{\bf Example 4.} Consider a quadratic ($k=2$) matrix polynomial
$P(\la)=A_2 \la^2+A_1 \la+A_0$, where

\begin{equation}\label{}
\begin{array}{ccc}
A_2=\left[
\begin{array}{cc}
  0        &      0 \\
      -2     &         2  \\
\end{array} \right],
A_1=\left[
\begin{array}{cc}
 1/7    &       -1/7 \\
      27/7     &     -41/7  \\
\end{array} \right],
A_0=\left[
\begin{array}{cc}
 2/7       &    -2/7 \\
     174/7     &    -146/7  \\
\end{array} \right]
\end{array}
\end{equation}

are matrices of order 2 ($n=2)$.

We have that $([X,Y],T\oplus Z)$, with
$$X=\left[
\begin{array}{cc}
1 & 1  \\
1 & 8  \\
\end{array} \right],
T=\left[
\begin{array}{cc}
2 & 0  \\
0 & -2  \\
\end{array} \right],
Y=\left[
\begin{array}{cc}
1 & 1  \\
1 & 2  \\
\end{array} \right],
Z=\left[
\begin{array}{cc}
0 & 1  \\
0 & 0  \\
\end{array} \right]
$$

 is a maximal standard pair for
$P(\lambda)$, with $Z$ nilpotent. $P(\la)$ has only one solvent

\begin{equation}\label{}
S_{1}= XTX^{-1}=\left[
\begin{array}{cc}
  18/7         &  -4/7 \\
      32/7      &    -18/7      \\
\end{array} \right].
\end{equation}

We can compute the cosolvents performing an spectral inversion,
$(Y',Z')$, with
$$Y'=\left[
\begin{array}{cccc}
1 & 1 & 1 & 1 \\
1 & 8 & 1 & 2  \\
\end{array} \right],
Z'=\left[
\begin{array}{cccc}
1/2 & 0 & 0 & 0 \\
0 & -1/2  & 0 & 0 \\
0 & 0 & 0 & 1 \\
0 & 0 & 0 & 0 \\
\end{array} \right].
$$
Therefore, with the leading vectors of the Jordan chains in $Y'$
we can construct 3 cosolvents
$$
\tilde S_{1}= \left[
\begin{array}{cc}
1 & 1  \\
1    & 8  \\
\end{array} \right]\left[
\begin{array}{cc}
1/2 & 0  \\
0    & -1/2  \\
\end{array} \right]\left[
\begin{array}{cc}
1 & 1  \\
1    & 8  \\
\end{array} \right]^{-1}=\left[
\begin{array}{cc}
   9/14       &   -1/7 \\
       8/7   &        -9/14  \\
\end{array} \right],$$

$$ \tilde S_{2}=\left[
\begin{array}{cc}
 1       &       1 \\
       8  &            1  \\
\end{array} \right]
 \left[
\begin{array}{cc}
-1/2 & 0  \\
0 & 0  \\
\end{array} \right]
\left[
\begin{array}{cc}
1       &       1 \\
       8  &            1  \\
\end{array} \right]^{-1}=
 \left[
\begin{array}{cc}
 1/14      &    -1/14 \\
       4/7  &         -4/7 \\
\end{array} \right]$$
and

$$\tilde S_{3}=\left[
\begin{array}{cc}
1 & 1  \\
1    & 2  \\
\end{array} \right]
 \left[
\begin{array}{cc}
0 & 1  \\
0 & 0  \\
\end{array} \right]
\left[
\begin{array}{cc}
1 & 1  \\
1   & 2  \\
\end{array} \right]^{-1}=
 \left[
\begin{array}{cc}
  -1       &       1 \\
      -1    &          1 \\
\end{array} \right],
$$
the cosolvent $\tilde S_3$ is nilpotent and so it is an "infinite
solvent" of $P(\la)$).

\section{Bisolvents and regular right factors}

Observe that a regular matrix polynomial may have regular first
order factors without admitting any solvent or cosolvent.
 We now extend solvent theory so as to yield regular right factors in some of such cases.

\begin{quote}\begin{definition}
\quad (i) A bisolvent for the $m\times m$ matrix $P(\la)$ is a
pair $(S_1,S_2)$ of commuting $m\times m$ matrices which satisfy
\beq\label{sol-sys} \quad(ii)\jump
\sum_{i=0}^kA_iS_1^iS_2^{k-i}=0.\eeq

\quad (iii) We call the bisolvent separable if for some $m\times
m$ idempotent $\Pi$ we have
  \beq\label{sep-sys} \jump S_1=\Pi S_1\Pi+(I-\Pi),\jump S_2=\Pi+(I-\Pi)S_2(I-\Pi).\eeq

\quad $(S_1,S_2,\Pi)$ is a solvent if $\Pi=I$ and cosolvent if
$\Pi=0.$

\end{definition}\end{quote}

Observe that, any pair $(S_1,S_2)$ such that $S_1S_2=0$ and both
$S_1$ and $S_2$ are nilpotent with an index $(k-1)$ satisfies item
(ii) of this definition, in special $(S_1=0,S_2=0)$. We are
interested in bisolvents that fulfill the item (iii), we will see
that in this case $S(\la)=\la S_2-S_1$ has spectral information of
$P(\la)$ and are right divisors.

 Observe that, in a separable bisolvent, $S_1$ and
$S_2$ normally fail to be a solvent or a cosolvent. Nevertheless,
solvents and cosolvents can be viewed as special cases: $\Pi=I$
and $S_2=I$ for solvents, $\Pi=0$ and $S_1=I$ for cosolvents.

As we shall see, the commutativity assumption does not limit the
scope of bisolvent theory. The separability hypothesis implies
that the pencil $\la S_1+S_2$ is commutative and regular.

\noindent The special case $P(\la)=\la I_n-A$ is instructive. Here, the roots
of $P(\la)$ coincide with the eigenvalues of $A$, and the corresponding Jordan basis of "eigenvectors"
for $P(\la)$ are the actual (nominal or
generalized) eigenvectors of $A$. Let $J$ be the
Jordan form for $A$ and let $X$ be the similarity matrix between them: $AX=XJ.$ Then
$(X,J)$ is a standard pair for $P(\la)$. Another one (the companion standard pair)
is $(I,A)$. The various standard pairs form a
single similarity orbit through the gauge transformation $(X,T)\ra(XG,G\inv TG).$

So, in this special case, $P(\la)$ admits many standard pairs but
only a single solvent, viz. $A$. The proceeding in Theorem
\ref{monic-solv-thm} (i) for producing the solvent from a standard
pair reduces to the equation $A=XTX\inv.$

The following procedure defines potential separable bisolvents for
$P(\la),$ and at the same time, potential right divisors.

Let $([X,Y],T\oplus Z)$ be any sufficiently large standard pair
for $P(\la)$. Choose an $m$-dimensional $T\oplus Z$-invariant
subspace $V_1\oplus V_2$ and let $(\tilde X,\tilde Y,\tilde
T,\tilde Z)$ be the corresponding restriction:
 $$\tilde X=X_{_{|V_1}},\quad \tilde Y=Y_{_{|V_2}},\quad \tilde T=P_{V_1}T_{_{|V_1}},\quad \tilde Z=P_{V_2}Z_{_{|V_2}}.$$
By specifying bases in $V_1$ and $V_2$ we may consider both
$Q_0=[\tilde X,\tilde Y]$ and $\tilde T\oplus \tilde Z$ as
$m\times m$ matrices. In the same basis, $I\oplus 0$ represents an
idempotent, i.e. the projection onto $V_1$ and along $V_2$. Now,
assuming $Q_0$ non-singular, define
 $$S_1=\tilde X(\tilde T_1\oplus I)\tilde X\inv, \quad
 S_2=\tilde X(I\oplus\tilde T_2)\tilde X\inv,\quad
 \Pi=\tilde X(I\oplus 0)\tilde X\inv.$$









This shows that indeed we have a separable bisolvent, and a right
factor, which can be enunciate as follows.

\begin{quote} \begin{lemma}\label{reg-bis} As long as $Q_0$ is
invertible, $(S_1,S_2)$ is a separable bisolvent and
$S(\la)=S_1\la-S_2$ is a regular first order right factor.
\end{lemma}\end{quote}

 {\bf Proof.}
$$\sum_{i=0}^kA_iS_1^iS_2^{k-i}=[\sum_{i=0}^kA_i\tilde X\tilde T^i
| \sum_{i=0}^kA_i\tilde Y\tilde
Z^{k-i}]Q_0\inv=[\sum_{i=0}^kA_iXT^i_{|V_1} |
\sum_{i=0}^kA_iYZ^{k-i}]_{|V_2}]Q_0\inv,$$ and the RHS vanishes
since $([X,Y],T\oplus Z)$ is a standard pair. Thus, $(S_1,S_2)$ is
a bisolvent, and clearly with $\Pi=I_p \oplus 0_{m-p}$, for a
suitable $p$, it is separable. \qed

The following results show that separable bisolvent theory is a
full-fledged generalization of solvent/cosolvent theories, and an
adequate tool in handling the regular case:

\begin{quote}\begin{theorem}\label{reg-bisolv-thm}
\quad (i)\quad  All the separable bisolvents $S$ of $P(\la)$ can
be constructed from a "maximal controllable bistandard pair"
$([X,Y],T\oplus Z)$ by restriction to an $m$-dimensional $T\oplus
Z$-invariant subspace. Namely,
 $$S_1\la-S_2=[\tilde X,\tilde Y][(I\la-\tilde T)\oplus(I-\la\tilde Z)][\tilde
 X,\tilde Y]\inv$$
where $([\tilde X,\tilde Y],\tilde T\oplus\tilde Z)$ are the
restrictions of $([X,Y],T\oplus Z)$ to the subspace.

\quad (ii) \quad In particular, the KW spectrum of $S_1\la-S_2$ is
subordinate to that of $(I\la-T)\oplus(I-Z\la)$.

\quad (iii) \quad $\la A+B$ is a right regular factor iff
$(S_1,S_2)=(XA,XB)$ is a separable bisolvent for some $X\in GL_m.$
In other words, separable bisolvents describe all the right
factors up to left equivalence!
\end{theorem}\end{quote}


 Item (iii) is central in establishing the natural
status of bisolvent theory in general factorization procedures. It
also implies a simplification in the calculation of the left
factor $Q(\la)=P(\la)(S_1\la-S_2)\inv$ for a given right factor.

Important conclusions can be derived from previous results,
suppose that $P(\la)$ is a regular matrix polynomial, let
$\mathcal{B}$, $\mathcal{BS}$ and $\mathcal{F}$ be the respective
sets of bisolvents, separable bisolvents and the regular first
order right factors, so $\mathcal{BS} \varsubsetneq \mathcal{B}$,
$\mathcal{BS} \varsubsetneq \mathcal{F}$ and moreover
$\mathcal{BS} = \mathcal{B} \cap \mathcal{F}$.

{\bf Example 5.} We demonstrate a subtlety which arises due to
Theorem \ref{reg-bisolv-thm} (ii) (spectral inversion). Consider
the matrix polynomial
\beq\label{pol1}P(\la)=\la\oplus(\la-1)\oplus
1=\la{\tiny\BpM 1&0&0\\ 0&1&0\\ 0&0&0\EpM}+{\tiny\BpM0&0&0\\
0&-1&0\\ 0&0&1\EpM},\eeq and with it two non-similar standard
pairs:
 \beq\label{bis1}X={\tiny\BpM 1&0\\ 0&1\\ 0&0\EpM},
 \quad Y={\tiny\BpM 0\\ 0\\ 1\EpM},
 \quad T=0\oplus 1, \quad Z=0\eeq
 \beq\label{bis1}\hat X={\tiny\BpM 1\\ 0\\ 0\EpM}, \quad \hat
 Y={\tiny\BpM 0&0\\ 1&0\\ 0&1\EpM},
 \quad \hat T=0, \quad \hat Z=1\oplus 0\eeq
related via a simple spectral inversion of the eigenvalue $\la=1.$
$P(\la)$ admits the bisolvent,
 $S_1=0\oplus 1\oplus 1,$ $S_2=1\oplus 1\oplus 0,$ which is separable w.r.t. two distinct idempotents,
  $\Pi=1\oplus 0\oplus 0$ and $\hat \Pi=1\oplus 1\oplus 0.$
It can be easily verified that $(S_1,S_2,\Pi)$ (resp.
$(S_1,S_2,\hat \Pi))$ is a restriction of  $(X,Y,T,Z)$ (resp.
$(\hat X,\hat Y,\hat T,\hat Z))$. In other words, spectral
inversion does not affect the bisolvent but may affect its
separating idempotent. \qed

The existence of a non-nilpotent solvent $S$ always leads a
separable bisolvent with an spectral inversion, with $\la S_1-
S_2= (\la I- S)Q$, with an appropriated  $Q$, the same can be
concluded with a cosolvent. Furthermore we have.

\begin{theorem}\label{solv-bisol}\quad (i) Let $(S_1,S_2)$ be a separable bisolvent if $S_1$ is
nonsingular,
 then $S=S_1^{-1}S_2$ is a solvent and if $S_2$ is nonsingular, then
$\tilde S=S_2^{-1}S_1$ is a cosolvent.

\quad (ii) In particular, if $(P(\la))$ has no zero nor infinite
eigenvalues, then $(S_1,S_2)$ is a separable bisolvent if and only
if $S=S_1^{-1}S_2$ is a solvent and if and only if $\tilde
S=S_2^{-1}S_1$ is a cosolvent.
\end{theorem} \end{quote}

In one word, in Item (ii) if $P(\la)$ has no solvents it will have
no cosolvents nor separable bisolvents, anyway this condition is
very restrictive, so we can conclude that the importance of
bisolvents appears when a zero eigenvalue or (and) an infinite
eigenvalue exist. We can suppose without loss of generality that
in Theorem \ref{reg-bisolv-thm} $Z$ is nilpotent, then it will be
clear that we can have separable bisolvents and therefore right
regular first order factors, with no solvents or (and) cosolvents
at all.

{\bf Example 6.} Consider a quadratic ($k=2$) matrix polynomial
$P(\la)=A_2 \la^2+A_1 \la+A_0$, where

\begin{equation}\label{}
\begin{array}{ccc}
A_2=\left[
\begin{array}{ccc}
 -33/10       &    3/2    &        8/5 \\
     -16/5     &      1   &           7/5 \\
      19/10     &   -1/2  &         -4/5    \\
\end{array} \right],
A_1=\left[
\begin{array}{ccc}

    -788/215    &     11/43   &      281/215 \\
     931/215     &   -22/43   &     -347/215  \\
    1097/215      & -134/43   &     -589/215   \\
\end{array} \right],
\end{array}
\end{equation}

\begin{equation}
A_0=\left[
\begin{array}{ccc}

     -60/43     &     76/43    &      74/43 \\
      30/43      &   -81/43    &       6/43  \\
     -36/43       &   80/43    &      10/43   \\
\end{array} \right]
\end{equation}
are matrices of order 2 ($n=3$).

We have that $([X,Y],T\oplus Z)$, with
$$ X=\left[
\begin{array}{ccc}
   2   &  2  &   5 \\
     4  &   4 &    2 \\
     1   &  1  &   2  \\
\end{array} \right],
T=\left[
\begin{array}{ccc}
       3      &        0    &          0 \\
       0       &       2     &         0 \\
       0        &      0      &        0 \\
\end{array} \right],
Y=\left[
\begin{array}{ccc}
    1   &  1  &   3 \\
    -1   & -1  &   1 \\
     3    & 3   &  4  \\
\end{array} \right],
Z=\left[
\begin{array}{ccc}
 0      &        1       &       0 \\
       0       &       0  &            1 \\
       0        &      0   &           0 \\
\end{array} \right]
$$

 is a maximal standard pair for
$P(\lambda)$, with $Z$ nilpotent. We can not construct any
nonsingular matrix of size $n$ with the leading vectors of $X$ nor
$Z$, therefore $P(\la)$ has no solvents nor cosolvents. Besides
that, with the standard pair above we get two separable
bisolvents: ($S_1,S_2$) and ($S'_1,S_2$), where
$S_1=U_1J_1U_1^{-1}$, $S_2=U_1J_2U_1^{-1}$ and
$S'_1=U_1J_3U_1^{-1}$, in which
 $$ U_1=\left[
\begin{array}{ccc}
 2      &    5    &       1 \\
     4    &   2   &        -1 \\
      1     &   2  &       3   \\
\end{array} \right],
J_1=\left[
\begin{array}{ccc}
    3      &        0    &          0 \\
       0       &       0     &         0 \\
       0        &      0      &        1 \\
\end{array} \right],
J_2=\left[
\begin{array}{ccc}
     1      &        0    &          0 \\
       0       &      1     &         0 \\
       0        &      0      &        0 \\
\end{array} \right],
J_3=\left[
\begin{array}{ccc}
    2      &        0    &          0 \\
       0       &       0     &         0 \\
       0        &      0      &        1 \\
\end{array} \right].
$$

\section{Bisolvents: the additive formalism}

Separable bisolvents can also be described using an equivalent additive formalism.

\begin{quote}\begin{definition}\label{def-add}
An additive (separable) bisolvent is a triple of $m\times m$ matrices
$(P_1,P_2,\Pi)$ with $\Pi$ idempotent, which satisfies
  \beq\label{add-sol}(i)\jump \sum_{i=0}^kA_i(P_1^i+P_2^{k-i})=0\eeq
 $$(ii)\jump P_1=\Pi T_1\Pi,\jump P_2=(I-\Pi)T_2(I-\Pi).$$
 \end{definition}\end{quote}

The two separable theories, multiplicative and additive, are completely equivalent.
Indeed, given $P_1,P_2,\Pi$ we define $S_1,S_2$ via
 $$S_1=P_1+(I-\Pi),\jump S_2=P_2+\Pi.$$
Conversely, given $S_1,S_2,\Pi$ we define
 $$P_1=\Pi S_1\Pi,\jump P_2=(I-\Pi)S_2(I-\Pi).$$
In both cases we have $P_1^i+P_2^j=S_1^iS_2^j$ holds for all $i,j\geq 0.$

Both bisolvent pairs are commutative. In the additive formalism, in fact, $P_1P_2=P_2P_1=0.$

\section{Conclusions}

The new bisolvent formulation can be an important contribution to
the factorization theory on regular matrix polynomials permitting
the description, up to equivalence, of all right regular first
order factors in contrast with solvents/cosolvents which restrict
it to monic/comonic factors.

\end{document}